\documentclass[11pt]{article}

\usepackage{amssymb}
\usepackage{amsfonts}
\usepackage{eucal}
\usepackage{soul}
\usepackage{amsmath}
\newtheorem{theorem}{Theorem}[section]
\newtheorem{proposition}[theorem]{Proposition}

\newtheorem{corollary}[theorem]{Corollary}

\newtheorem{lemma}[theorem]{Lemma}
\newtheorem{remark}[theorem]{Remark}

\newenvironment{proof}[1][Proof]{\textbf{#1.} }{\ \rule{0.5em}{0.5em}}
\input pictex.sty

\usepackage{cancel}
\usepackage{xcolor}

\usepackage{amsmath,color,ulem}
\newcommand{\stkout}[1]{\ifmmode\text{\sout{\ensuremath{#1}}}\else\sout{#1}\fi}

\begin{document}

 \title{Strong laws of large numbers for lightly trimmed sums of generalized Oppenheim expansions}
\author{Rita Giuliano \footnote{Dipartimento di
		Matematica, Universit\`a di Pisa, Largo Bruno
		Pontecorvo 5, I-56127 Pisa, Italy (email: rita.giuliano@unipi.it)}~~and Milto Hadjikyriakou\footnote{School of Sciences, University of Central Lancashire, Cyprus campus, 12-14 University Avenue, Pyla, 7080 Larnaka, Cyprus (email:
		mhadjikyriakou@uclan.ac.uk). Part of this work was conducted while the author was a visiting scholar at the University of Cyprus.   }		
		\footnote{The authors wish to thank I. Berkes for having drawn   to their attention the subject of the present work. } } 	 

\maketitle

\begin{abstract}\noindent In the framework of generalized Oppenheim expansions we prove strong law of large numbers for lightly trimmed sums. In the first part of this work we identify a particular class of expansions for which we provide a convergence result assuming that only the largest summand is deleted from the sum; this result  generalizes a  strong law  recently proven for the L\"uroth case. In the second part we drop any assumptions concerning the structure of the Oppenheim expansions and we prove  a result concerning trimmed sums when  at least  two  summands are trimmed; then we derive a corollary for the case  in which only the largest summand is deleted from the sum.


\end{abstract}

	

\section{Introduction} 
The framework of this work has been introduced and generalized in \cite{G2018} and \cite{GH2020} respectively and is described as follows: let $(B_n)_{n\geq  1}$ be a sequence of integer-valued random  variables defined on $(\Omega, \mathcal{A}, P)$, where $\Omega =[0,1]$, $\mathcal{A}$ is the $\sigma$-algebra of the Borel subsets of $[0,1]$ and $P$ is the Lebesgue measure on $[0,1]$. Let $\{F_n, n\geq 1\}$ be a sequence of probability distribution functions   with  $F_n(0)=0$, $F_n(1)=1$ $\forall n$ and moreover let $\varphi_n:\mathbb{N}^*\to \mathbb{R}^+$ be a sequence of functions. Furthermore, let $(q_n)_{n\geq  1}$ with $q_n=q_n(h_1, \dots, h_n)$ be a sequence of nonnegative numbers (i.e. possibly depending on the $n$ integers $h_1, \dots, h_n$) such that, for $h_1 \geq  1$ and $h_j\geq  \varphi_{j-1}(h_{j-1})$, $j=2, \dots, n$ we have
\[
P\big(B_{n+1}=h_{n+1}|B_{n}=h_{n}, \dots, B_{1}=h_{1}\big)= F_n(\beta_n)-F_n(\alpha_n),
\]
where 
\[
\alpha_n=\delta_n(h_n, h_{n+1}+1, q_n)  ,\quad \beta_n=\delta_n(h_n, h_{n+1}, q_n)\quad\mbox{with}\quad\delta_j(h,k, q) = \frac{ \varphi_j (h )(1+q )}{k+\varphi_j (h ) q  }.
\]
Let $Q_n= q_n(B_1, \dots, B_n)$ and define
\begin{equation} 
\label{Rdef}R_{n}= \frac{ B_{n+1}+\varphi_n(B_n) Q_n}{\varphi_n(B_n)(1+Q_n) }= \frac{1}{\delta_n(B_n, B_{n+1}, Q_n)}.
\end{equation}

\noindent In \cite{GH2020} (see Lemma 3 there) it has been proven that for any integer $n$ and for $x\geq 1$, 
\begin{equation}\label{A4}  
P(R_n>x)\leq F_n\left(\frac{1}{x}\right)
\end{equation}
which implies that if $U_n$ is a  random variable with distribution $F_n$ and $Y_n = \displaystyle \frac{1}{U_n}$ for every integer $n$, then
\begin{equation}
\label{stdom}P(R_n>x)\leq P(Y_n>x),\quad \forall x\geq 1.
\end{equation}
i.e., the sequence $(R_n)_{n\geq 1}$ is \textit{stochastically dominated} by the sequence $(Y_n)_{n\geq 1}$. 

\medskip

Since the random variables $(R_n)_{n\geq 1}$ do not have finite expectations, a traditional strong law for the quantity $\frac{1}{a_n}\sum_{i=1}^{n}R_i$ cannot be proven. However, in \cite{G2018}, under some conditions for the involved distributions, the convergence in probability of $\dfrac{1}{n\log n}\sum_{i=1}^{n}R_i$ is established. This result, raises the question whether a strong law of large numbers can be proven, after deleting finitely many of the largest summands from the partial sums. Particularly,   let $r$ be a fixed integer. We are interested in studying the almost sure convergence of 
\[
\dfrac{{}^{(r)}S_n}{n\log n},
\]
where  ${}^{(r)}S_n = \sum_{i=1}^nR_i - \sum_{k=1}^{r}M_n^{(k)}$, $M_n^{(r)}$ denoting the $r$-th maximum of $R_1,\ldots,R_n$ (in decreasing order i.e. $M_n^{(1)}$ denotes the maximum). In the literature,  the sequence  $\big({}^{(r)}S_n\big)_{n\geq 1}$ is known as the \textit{lightly trimmed sum process}.  Note that in the case where $r$ is substituted by a sequence $(r_n)_{n\geq 1}$ such that $r_n \to \infty$ and $r_n/n\to 0$ as $n\to \infty$ we have the so-called \textit{moderate trimming} while in the case where $r_n/n\to c\in (0,1)$ the resulting sequence is said to be \textit{heavily trimmed}. For more details we refer the interested reader to \cite{BHS2012} (and references therein). Convergence results for moderately trimmed sums of Oppenheim expansions can be found \cite{GH2023}.  

The structure of the paper is as follows: in Section \ref{sec: specialcase} we introduce a special class of Oppenheim expansions for which we establish a strong law for its lightly trimmed sum processes for the case $r=1$; this  class enjoys some particular  characteristics   that allow  us  to apply known results for independent random variables, and our Theorem \ref{Th1} generalizes a   recent  result  obtained in  \cite{AA} for the L\"uroth expansion. In Section \ref{sec: generalcase} we address the general case (i.e. we do not impose any condition on the sequence of  expansions taken into account), and we have been able to prove   an asymptotic result for $r\geq 2$. This almost sure convergence result becomes instrumental for proving another asymptotic law for $r=1$ which, though being weaker than the one presented in   Section \ref{sec: specialcase},  is  in the same direction. It is worth noting that for this latter result any assumptions for the structure of the Oppenheim expansion are dropped and we impose more relaxed conditions for the involved distribution functions.

\section{A strong law for a class of generalized Oppenheim expansions  and  $r=1$}
\label{sec: specialcase}

\noindent For the random variables $(R_n)_{n\geq 1}$ defined above, the following two relations were proven (see Lemma 2 relation (5) and Lemma 3 respectively in \cite{GH2020}): for $x, y \geq 1$ and $m<n$, 
$$P(R_m > x)= E\Big[F_m\Big(\frac{\phi_m(B_m)(1+Y_m)}{\lceil x \phi_m(B_m)+(x-1)Y_m\phi_m(B_m)\rceil +Y_m\phi_m(B_m) }\Big)\Big],$$
$$P(R_m > x, R_{n}> y)=E\Big[I(R_m > x)F_n\Big(\frac{\phi_n(B_n)(1+Y_n)}{\lceil y \phi_n(B_n)+(y-1)Y_n\phi_n(B_n)\rceil +Y_n\phi_n(B_n) }\Big)\Big],$$
where $\lceil x\rceil$ denotes the least integer greater than or equal to $x$.
Then, the following proposition is obvious.

\begin{proposition} \sl\label{propB} 
	For the random variables $(R_n)_{n\geq 1}$, the following results hold true:
	\begin{enumerate}
		\item[(a)] Assume that $x\geq 1$ is such that $x \phi_m(B_m)+(x-1)Y_m\phi_m(B_m)$ is an integer for any integer $m$. Then
		$$P(R_m > x)= F_m\Big(\frac{1}{x}\Big), \, \forall x\geq 1.$$
		
		\item [	(b) ] Assume in addition that $y\geq 1$ is  such that $y \phi_n(B_n)+(y-1)Y_n\phi_n(B_n)$ is an integer for any integer $n$. Then
		$$P(R_m > x, R_{n}> y)=F_m\Big(\frac{1}{x}\Big)F_n\Big(\frac{1}{y}\Big),  \, \forall x,y\geq 1.$$
	\end{enumerate}

\end{proposition}

\begin{proposition} \sl\label{propA}Consider the random variables $(R_n)_{n\geq 1}$ and assume that $x_i\geq 1, \, \forall i=1,2,\ldots,n$. Then,
	\begin{align*}&
	P(R_1>x_1,\ldots, R_n>x_n)\\& = E\left[ F_n\left( \frac{\phi_n(B_n)(1+Y_n)}{\lceil x_n\phi_n(B_n)+(x_n-1)Y_n\phi_n(B_n)\rceil +\phi_n(B_n)Y_n}\right)I(R_1>x_1,\ldots, R_{n-1}>x_{n-1})\right].
	\end{align*}
\end{proposition}
\begin{proof}
	The proof follows by applying similar steps as in the case of $n=2$ (Lemmas 2 and 3 in \cite{GH2020}).
\end{proof}

\begin{proposition}
	\label{propExt} \sl  Consider the random variables $(R_n)_{n\geq 1}$. Then, for every integer $n$ and for every finite set of numbers $x_i\geq 1, \, \forall i=1,2,\ldots,n$ such that  $x_k\phi_k(B_k)+(x_k-1)Y_k\phi_k(B_k)$ is an integer for every $k = 1, 2, \ldots, n$, we have 
	\[
	P(R_1>x_1,\ldots, R_n>x_n) = F_1\left(\frac{1}{x_1}\right)\ldots F_n\left(\frac{1}{x_n}\right).
	\]
\end{proposition}
\begin{proof}
	The result follows by induction.  The case  $n=2$ is discussed in Proposition \ref{propB}. Assume that the statement is true for $n-1$. Then by Proposition \ref{propA} we can write
	\begin{align*}
	&P(R_1>x_1,\ldots, R_n>x_n) \\ &= E\left[ F_n\left( \frac{\phi_n(B_n)(1+Y_n)}{\lceil x_n\phi_n(B_n)+(x_n-1)Y_n \phi_n(B_n)\rceil +\phi_n(B_n)Y_n}\right)I(R_1>x_1,\ldots, R_{n-1}>x_{n-1})\right]\\
	& = F_n\left(\frac{1}{x_n}\right)P(R_1>x_1,\ldots, R_{n-1}>x_{n-1})
	\end{align*}
	which  leads to the conclusion, by the induction hypothesis.
\end{proof} 

\medskip

The proposition that follows will be a ``key" result for obtaining the convergence theorem of this section: by using a particular class of Oppenheim expansions we define a sequence of discrete random variables that is proven to consist of independent random variables the densities of which can be easily calculated.

\medskip
\noindent 
We start with the definition of a sequence which plays an important role for the results of this section.

\bigskip
\noindent  Call {\it good} a strictly increasing sequence  $\Lambda = (\lambda_j)_{j \in \mathbb{N}}$    tending to $+ \infty$  with $\lambda_j\geq 1$ for every $j\geq 1$ and $\lambda_0 =0$. For $u \in [1, + \infty)$ let $j_u$ be the only integer such that $\lambda_{j_u-1}< u \leq \lambda_{j_u}$  (i.e. $\lambda_{j_u}$ is the minimum element in $\Lambda$ larger than or equal to $u$).
  
\begin{proposition}\sl\label{indepTn} 
Consider the random variables $(R_n)_{n\geq 1}$ and assume that there exists a good sequence $\Lambda$  such that
  for every $x \in \Lambda$ and for every $n$,  $x \phi_n(B_n)+(x -1)Y_n\phi_n(B_n)$ is an integer.For every $n$, denote
  \begin{equation}\label{Tn}
  T_n = \lambda_{j_{R_n}}.
  \end{equation}
Then $T_n$ takes values in $\Lambda$, and the sequence $(T_n)_{n \geq 1}$ consists of independent random variables. Moreover the discrete density of $T_n$  is  given by the formula $$F_n\Big(\frac{1}{\lambda_{s-1}}\Big)-F_n\Big(\frac{1}{\lambda_s}\Big),\quad s \in \mathbb{N}^*.$$
 \end{proposition} 
\begin{proof}
	Observe that, for any integer $n$,  the relation $\lambda_{j_r} >\lambda_n \Leftrightarrow r >\lambda_n $. Thus, for every finite set of integers $\{i_1, \dots, i_k\}$ and for every finite set of integers $n_{i_1}, \dots, n_{i_k}$ we have
	$$P(T_{i_1}> \lambda_{n_{i_1}}, \dots, T_{i_k}> \lambda_{n_{i_k}})=P(R_{i_1}>  \lambda_{n_{i_1}} , \dots, R_{i_k}>  \lambda_{n_{i_k}} )= F_{i_1}\Big(\frac{1}{ \lambda_{n_{i_1}}}\Big)\cdot \dots \cdot F_{i_k}\Big(\frac{1}{ \lambda_{n_{i_k}}}\Big),$$
	which proves the independence of the random variables $\{T_n\}_{n}$. For the density, note that for every integer $     s \in \mathbb{N}^*$, we have
	$$P(T_n = \lambda_s)= P(T_n >\lambda_{s-1})- P(T_n > \lambda_s)= F_n\Big(\frac{1}{\lambda_{s-1}}\Big)-F_n\Big(\frac{1}{ \lambda_s}\Big).$$	  
\end{proof}

   


\begin{remark}
	The result above is a generalization of Theorem 3 in Galambos in \cite{G1974},  in which $y_n = 0$, $\Lambda= \mathbb{N}$ and $F_n(x) = F(x) = x$. 
\end{remark}

\begin{remark}\label{remark1}
	It is important to identify functions $\phi_n$ for which the conditions imposed in Proposition \ref{indepTn} are satisfied. First, recall that the notation $y_n$ stands for the sequence of nonnegative numbers such that $y_n(B_1, \dots, B_n)= Y_n$. 	 As a first example, consider positive integers $a_1, \dots, a_p$ and assume that   
	\[
	\phi_{kp +j-1} = 1/a_j , \quad \mbox{for  }\quad k \in \mathbb{N}, \quad j=1, \dots, p \\
	 \]
	 Define $\kappa = L.C.M. (a_1, \dots, a_p)$ and $\Lambda = ( \kappa n )_{n\geq 1}$ and
	assume that $y_n \equiv c_n$ where $(c_n)_{n \geq 1}$ is a sequence of positive numbers chosen from the  set $\Lambda$. Then, for any $x\in \Lambda$, $$x\phi_n(B_n) + (x-1)Y_n\phi_n(B_n)$$ is an integer. 
	
	\medskip
\noindent	
Moreover, the conditions of the proposition are satisfied if $\Lambda = \mathbb{N^*}$, $y_n \equiv 0$ and $\phi_n(h) = \sum_{k=1}^m h^k$ for some integer $m$. Note that for $m= 1$ we get the corresponding $\phi$ function for the Engel series while for $m=2$ we have the Sylvester expansion. If $\phi_n (h) = \sum_{k=0}^m h^k$ the case of $m = 0$ covers the L\"{u}roth case (see \cite{G2018} for details). 
\end{remark}

\medskip
\noindent	
Before stating and proving the main result of this section, we present, without a proof, the result below which is part of Theorem 1 in \cite{M} and it is instrumental for the proof of the convergence result we are interested in.

\begin{theorem} \sl \label{Mori} Let $(X_n)_{n \geq 1}$ be a sequence of i.i.d. random variables and denote $^{(r)}S_n$   the $n$-th sample sum  with the first $r$ largest terms removed.  Let $A$ be an absolutely continuous increasing function  defined on $[0,+\infty)$,  with $A(0)=0$ and satisfying
	\begin{itemize}
		\item[(i)] $\frac{A(x)}{x^{\frac{1}{\alpha}}}$ is non decreasing for some $\alpha \in (0,2)$,
		\item[(ii)] $\sup_{x >0} \frac{A(2x)}{A(x)}< \infty$,	
	\end{itemize}
	and let $B$ be its inverse function. For every $s>0$, denote $$J_s =  \int_{1}^{\infty}[P(X_1>x)]^s{\rm d}B^s(x)  $$  
	and assume that  $J_{r+1}< +\infty $; then there exists a sequence $(c_n)_{n\in\mathbb{N}}$ of numbers such that
	$$\lim_{n \to \infty}\frac{ ^{(r)}S_n}{A(n)}- c_n =0, \qquad P-a.s.$$
	Moreover the constants $c_n$ can be chosen to be
	$$c_n = \frac{n}{A(n)}\int_{|x|\leq \tau A(n)} x {\rm d} \{  P(X_1 \leq x )\},$$
	where $\tau>0$ is an arbitrary constant. 
\end{theorem} 
 
\medskip
\noindent	
Next, we state and prove the main result of this section, which is a strong law for the trimmed sums of the special class of generalized Oppenheim expansions discussed above in the case where only the maximum term is excluded i.e. for $^{(1)}S_n$. This result covers the L\"{uroth} case studied in \cite{AA} but it can also be used to derive the respective convergence for the Engel and Sylvester series (see Remark \ref{remark1}). To the best of our knowledge, the results for the Engel and Sylvester expansions are new.

\begin{theorem} 
	\label{Th1}  \sl Let the assumptions of Proposition \ref{indepTn} hold and moreover assume  the following:
\begin{enumerate}
	\item [(i)]  \begin{equation}\label{as1F}
	\sup_n (\lambda_{n+1}- \lambda_n)= \ell< + \infty;
	\end{equation}
	
	\item[(ii)] $F_n \equiv F$ for all integers $n$ and there exists a constant $ \alpha >0$ such that 
	\begin{equation}
	\label{as3F} \lim_{t\to 0}  \frac{F (t)}{t}= \alpha.
	\end{equation}
\end{enumerate}
Then,
	\[
	\dfrac{S_n -M_n^{(1)}}{n\log n} \to \alpha \quad\mbox{a.s.}
	\]
	where $S_n = \displaystyle\sum_{i=1}^{n}R_i$ and $M_n^{(1)} = \max\{R_1\ldots, R_n\}$.
\end{theorem}
\begin{proof}
	Following the notation introduced earlier, let $T_n = \lambda_{j_{R_n}} $ and define $\tilde M_n^{(1)}  = \max\{T_1\ldots, T_n\}$. Then,
	\[
	T_n-\ell\leq R_n \leq T_n\qquad\mbox{and}\qquad \tilde M_n^{(1)} -\ell \leq M_n^{(1)}\leq \tilde M_n^{(1)}.
	\]
	Thus,
	\[
	\frac{\sum_{k=1}^{n} T_k - \tilde M_n^{(1)} -\ell n }{n\log n} \leq \dfrac{S_n -M_n^{(1)}}{n\log n} \leq \frac{\sum_{k=1}^{n} T_k - (\tilde M_n^{(1)} -\ell) }{n\log n}
	\]
	so it is sufficient to study the convergence of 
	\[
	\frac{\sum_{k=1}^{n} T_k - \tilde M_n^{(1)} }{n\log n}.
	\]
	By Proposition \ref{indepTn}, the sequence $(T_n)_{n\geq 1}$ defined in \eqref{Tn} consists of independent and identically distributed random variables and therefore Theorem \ref{Mori} can be employed; to this extent, since we are interested in the case $r=1$ and $A(x) = x \log x$,   first we have to check that 
	\[
	J_2 = \int_{1}^{\infty}[P(T_1>x)]^2{\rm d}B^2(x) <\infty 
	\]
	where $B(x)$ is the inverse of $A(x) = x\log x$. Note that ${\rm d}B^2(x) = 2B(x)B'(x){\rm d}x$ while ${P(T_1>x) = F\left(\frac{1}{x}\right)}$ (by Proposition \ref{indepTn}). Hence, due to \eqref{as3F}, we have that 
	\[
	J_2 = \int_{1}^{\infty}F^2\left(\frac{1}{x}\right)2B(x)B'(x){\rm d}x\leq C_1 + C_2\int_{1}^{\infty}\left(\frac{1}{x^2}\right)B(x)B'(x){\rm d}x.
	\] 
Now use the change of variables $B(x) = y$; since $x=A(y)$ and ${\rm d}y = B'(x){\rm d}x$ we have that 
	\[
	J_2 \leq C_1 + C_2\int_{B(1)}^{\infty} \frac{y}{A^2(y)}{\rm d}y = \frac{C_1 + C_2}{\log B(1)}<\infty.
	\]
	Hence, by Theorem \ref{Mori}, there is $c_n$ such that as $n\to \infty$
	\[
	\frac{\sum_{k=1}^{n} T_k - \tilde M_n^{(1)} }{n\log n}-c_n\to 0\qquad P-a.s., 
	\]
	where 
	\[
	c_n = \frac{n}{A(n)}\int_{1}^{A(n)}x{\rm d} \{  P(T_1 \leq x )\} = -\dfrac{1}{\log n}\int_{1}^{n\log n}x{\rm d} \{ P(T_1 > x )\} = -\dfrac{1}{\log n}\int_{1}^{n\log n}x{\rm d}  F\Big(\frac{1}{x}\Big). 
	\]	
	Using integration by parts we have that
	\[
	c_n = -\dfrac{1}{\log n}\left( n\log nF\Big(\frac{1}{n \log n}\Big)- 1 \right) +\frac{1}{\log n}\int_{1}^{n\log n}F\left(\dfrac{1}{x}\right){\rm d}x 
	\]	
	which can be equivalently written as
	\[
	c_n =  -\dfrac{1}{\log n}\left(   \dfrac{F \left(\dfrac{1}{n\log n}\right)}{\dfrac{1}{n\log n}}     \right) +\dfrac{  1}{\log n}+\frac{1}{\log n}\int_{1}^{n\log n}F\left(\dfrac{1}{x}\right){\rm d}x =I_1+I_2+I_3.
	\]
	Obviously $I_2 \to 0$ and by employing \eqref{as3F} we have that $I_1 \to 0$ for $n\to\infty$. For $I_3$, we start by observing that
	\[
	\int_{1}^{n\log n}F\left(\dfrac{1}{x}\right){\rm d}x = \int_{1}^{\frac{1}{n\log n}}-\dfrac{F(y)}{y^2}{\rm d}y = \int_{\frac{1}{n\log n}}^{1}\dfrac{F(y)}{y^2}{\rm d}y.
	\] 
	By \eqref{as3F}, for fixed $\epsilon>0$ let $\delta\in (0,1)$ be such that,  for $y \in (0,\delta),$
	\[
	\alpha-\epsilon\leq \dfrac{F(y)}{y}\leq\alpha+\epsilon
	\]
	and let $n_0$ be sufficiently large such that $\frac{1}{n\log n}<\delta$ for $\forall n\geq n_0$. Then
	\[
	\int_{\frac{1}{n\log n}}^{\delta}\dfrac{(\alpha-\epsilon)}{y}{\rm d}y<  \int_{\frac{1}{n\log n}}^{\delta}\dfrac{F(y)}{y^2}{\rm d}y<\int_{\frac{1}{n\log n}}^{\delta}\dfrac{(\alpha+\epsilon)}{y}{\rm d}y
	\]
	which leads to 
	\[
	(\alpha-\epsilon)\log\delta-(\alpha-\epsilon)\log\left(\frac{1}{n\log n}\right)<\int_{\frac{1}{n\log n}}^{\delta}\dfrac{F(y)}{y^2}{\rm d}y<	(\alpha+\epsilon)\log\delta-(\alpha+\epsilon)\log\left(\frac{1}{n\log n}\right).
	\]
Hence, due to the arbitraryness of $\epsilon$, 
	\begin{equation} \label{A5}
	\int_{\frac{1}{n\log n}}^{\delta}\dfrac{F(y)}{y^2}{\rm d}y\sim \alpha(\log n +\log(\log n)).
	\end{equation}
	Moreover, for $n\to \infty$
	\begin{equation} \label{A6}
	\dfrac{1}{\log n}\int_{\delta}^{1}\dfrac{F(y)}{y^2}{\rm d}y\to 0
	\end{equation}
	 Relations  \eqref{A5} and \eqref{A6} together give that $I_3 \to \alpha $ as $n\to\infty$. This concludes the proof.
\end{proof}

\section{A general strong law}
\label{sec: generalcase}

\bigskip
\noindent
In this section we prove a strong law for the lightly trimmed sums of any Oppenheim expansion. It is important to highlight that the result is obtained without any assumption   for the dependence structure of the random variables $R_n$. Moreover, the convergence result is obtained without assuming common law for the involved distributions.

\bigskip
\noindent
 The main result of this section is presented below. Observe that we are interested in studying the asymptotic behaviour of the trimmed sums of Oppenheim expansions if $r\geq 2$ i.e. $r$ of the maximum terms are removed from the sum. However, this result will lead to an asymptotic law for the case $r=1$ as well, which, although weaker than the one proven in the previous section, is obtained under less restrictive conditions.
 
 \medskip
\begin{theorem}\label{conv1}\sl Consider the random variables $(R_n)_{n\geq 1}$ and assume that for the involved distribution functions $(F_n)_{n\geq 1}$ the following condition is satisfied:
	\begin{equation}
	\label{A1} \sup_{n \geq 1 }\limsup_{x\to 0}  \frac{F_n(x)}{x}<\infty.
	\end{equation}Then, for every $p > 2$ and $r \geq 2$,
	\begin{equation}
	\label{claim1}\lim_{n \to \infty}\frac{^{(r)}S_n}{(n \log n)^p}=0, \qquad P-{\rm a.s.} 
	\end{equation}

\end{theorem}

\medskip
\noindent
 Throughout this  section we denote  $a_n = n\log n$; moreover, for a given positive increasing function $f$, and for every integer $n$, we let $m_n = \lfloor \log_2 n\rfloor $ where $\lfloor x \rfloor $ is the greatest integer less than or equal to $x$. Additionally, let $t_n = f(2^{m_n})$. 
 
\bigskip
\noindent
 Before proving the strong law we are interested in, we prove some preliminary results. 

\begin{lemma}\label{lemma1}\sl  Consider  the random variables  $(R_n)_{n\geq 1}$ and let the related distributions $(F_n)_{n\geq 1}$ satisfy assumption \eqref{A1}.  Moreover, assume that $f$ is a positive  increasing function such that, for some  integer $r>0$,
	\begin{equation} 
 \label{A2}	\sum_{m=1}^{\infty} \left( \frac{2^m}{f(2^m )}\right)^r=\sum_{m=1}^{\infty} \left( \frac{2^m}{t_{2^m }}\right)^r<\infty.
 	\end{equation}
	Then,
	\[
	P(M_n^{(r)}>t_n \, \, {\rm i.o.}) = 0.
	\]
\end{lemma}
\begin{proof}
	For any integer $j\geq 0$ we define the event $$A_j  = \{R_i >t_{2^j} \mbox{ for at least $r$ indices such that } i<2^{j+1}\}$$
	and, for any integer $n\geq 1$ 
	$$B_n =\{M_n^{(r)}>t_n\}.$$
	Let $j$ be fixed and note that, for every $n$  such that    $2^j \leq n < 2^{j+1}$, we have $B_n \subseteq A_j$; thus 
	$$\bigcup_{\{n: m_n =j\}} B_n \subseteq A_j,$$
	which implies that
	\begin{align*}&
	\{M_n^{(r)}> t_n \, \,{\rm i.o.}\}= \bigcap_s \bigcup_{n \geq  s}B_n \subseteq  \bigcap_s \bigcup_{\{n:m_n \geq  m_s\}}B_n = \bigcap_s \bigcup_{j \geq  m_s}\Big(\bigcup_{\{n: m_n =j\}} B_n\Big)\subseteq  \bigcap_s \bigcup_{j \geq m_s}A_j\\
	& =\bigcap_k \bigcup_{ j \geq  k }A_j= \{A_j \, \, {\rm i.o.}\}.
	\end{align*} 
	The first ``$\subseteq$"  holds true since for  $s \geq 1$ we have that  $\{n: n \geq s\}\subset\{n: m_n \geq m_s\}$, while the third equality is valid based on the below observation
	\begin{align*}& \bigcap_{s=1}^\infty \Big(\bigcup_{j \geq m_s}A_j\Big)= \bigcap_{k=0}^\infty \Big\{\bigcap_{s=2^k}^{2^{k+1}-1}\Big(\bigcup_{j \geq m_s}A_j\Big)\Big\}= \bigcap_{k=0}^\infty  \Big\{\bigcap_{s=2^k}^{2^{k+1}-1}\Big(\bigcup_{j \geq k}A_j\Big)\Big\}=\bigcap_{k=0}^\infty \Big(\bigcup_{j \geq k}A_j\Big).
	\end{align*} 
	 
	Now, for every integer $m$,
	\begin{align*}
	&P(A_m)\leq  \sum_{1\leq i_1<i_2<\ldots<i_r <2^{m+1}}P(R_{i_1}>t_{2^m }, \ldots,R_{i_r}>t_{2^m })\leq \sum_{1\leq i_1<i_2<\ldots<i_r<2^{m+1}} \prod_{j=1}^{r}F_{i_j}\left(\frac{1}{t_{2^m }}\right)\\
	&\leq C\sum_{1\leq i_1<i_2<\ldots<i_r<2^{m+1}}\dfrac{1}{ (t_{2^m })^r}\leq C {2^{m+1} \choose r}\dfrac{1}{ (t_{2^m })^r}\leq C\left(\frac{2^m}{ t_{2^m }}\right)^r,
	\end{align*}
	where the  second inequality is due to Proposition  \ref{propExt} and the third one   to condition \eqref{A1}.
	Thus,
	\[
 	\sum_{m=1}^{\infty}P(A_m) \leq C\sum_{m=1}^{\infty} \left(\frac{2^m}{t_{2^m }}\right)^r
	\]
	which is finite because of \eqref{A2}. The result follows by the Borel-Cantelli lemma.
\end{proof}

\medskip
\noindent
\begin{corollary}\label{corollary2}\sl  Under the assumptions of Lemma \ref{lemma1} and   for any $r\geq 2$,
$$\lim_{n \to \infty} \frac{M_n^{(r)}}{n \log n }=0, \qquad  P-{\rm a.s.}$$

\end{corollary}
\begin{proof}
  It is sufficient to prove that, for every $\varepsilon >0$, we have
\begin{equation}\label{A3}
P(M_n^{(r)}>\varepsilon n\log n \, \, {\rm i.o.}) = 0.
\end{equation}
 This  can be derived from the previous lemma   by choosing $f(x) = \varepsilon x \log x$  for $x\geq 1$ (recall that $f$ is assumed to be increasing), since   for this function we can easily  obtain $t_{n}\leq C \varepsilon n \log n $  and therefore 
	\[
	P(M_n^{(r)}>\varepsilon n\log n \, \, {\rm i.o.})\leq P(M_n^{(r)}>t_n \, \, {\rm i.o.}).
	\]
The conclusion follows by observing that 
$$\sum_{m=1}^{\infty} \left( \frac{2^m}{ t_{2^m }}\right)^r= C\sum_{m=1}^{\infty}\frac{1}{m^r} <\infty.$$
 \end{proof}

\medskip

\begin{lemma}\label{lemma5}\sl   Assume that the  conditions  of Lemma \ref{lemma1} are satisfied and additionally $f(x) > x$ ultimately. For every $m$, denote  by $N_m$  the number of indices $j<2^{m+1}$ for which $R_j> t_{2^m}$.
Then, for every integer  $s$ such that  $s\geq r$, 
$$P(N_{m_n}\geq s\, \,{\rm i.o.})=0.$$
\end{lemma}
\begin{proof} Observe that
$$  P(N_{m_n} \geq s\, \,{\rm i.o.})\leq P(M_n^{(s)}\geq t_{2^{m_n} }\, \,{\rm i.o.}),$$
and   we can apply the Borel-Cantelli lemma  because of Lemma \ref{lemma1}: in fact, ultimately,$$  \left( \frac{2^m}{f(2^m )}\right)^s \leq \left( \frac{2^m}{f(2^m )}\right)^r, $$
  whence
$$	\sum_{m=1}^{\infty}\left( \frac{2^m}{f(2^m )}\right)^s< \infty,$$
 due to \eqref{A2}.
\end{proof}

\medskip

\begin{lemma}\label{lemma4} \sl  
Let $p, q>0$ be fixed and consider the series 
\begin{equation}\label{series}
\sum_n\frac{1}{n^{2qp-2q +1}(\log n)^{2pq}}\sum_{j = 1}^n t^{2q}_j, 
\end{equation}
  where $t_n = 2^{m_n} (\log 2^{m_n})^\alpha$. Then,  
\begin{itemize}
\item[(i)]  if $p > 2+\frac{1}{2q}$, then \eqref{series} converges for every $\alpha>0$;
\item[(ii)] if $p =  2+\frac{1}{2q}$   then \eqref{series} converges for $\alpha<2$.
\end{itemize}
\end{lemma}
\begin{proof}
 Choose $f(x) = x \log^\alpha x$   for $x\geq 1$ and  $\alpha> 0$.    Then,   $t_n = f(2^{m_n}) =2^{m_n} (\log 2^{m_n})^\alpha$.   Observe that, 
$$\sum_{j = 1}^n t^{2q}_j\leq \sum_{k=0}^{m_n +1}\left(\sum_{j= 2^k }^{2^{k+1}-1 }f^{2q}(2^k)\right)=\sum_{k=0}^{m_n +1} 2^{k(2q+1)} (\log 2^k)^{2\alpha q}= C \sum_{k=0}^{m_n +1}2^{k(2q+1)}k ^{2\alpha q}.$$
By an application of Cesaro theorem, it is not difficult to see that 
$$\sum_{k=0}^{N} 2^{k(2q+1)}k ^{2\alpha q} \sim C\cdot  2^{N(2q+1)}N ^{2\alpha q}, \qquad N \to \infty.$$
Hence, ultimately
$$\sum_{j = 1}^n t^{2q}_j\leq C\cdot  2^{(2q+1)m_n }m_n ^{2\alpha q}\leq  C  \cdot 2^{{(2q+1)}\log_2 n }  (\log_2 n)^{2\alpha q}= C \cdot n^{2q+1} (\log n)^{2q\alpha},$$
and $$\frac{1}{n^{2qp-2q +1}(\log n)^{2qp}}\sum_{j = 1}^n t^2_j\leq  \frac{C}{n^{2pq-4q}(\log n)^{2q(p-\alpha)}}.$$
The claim follows by known results on Bertrand series.
\end{proof}

\medskip

\begin{corollary}\label{corollary1}\sl  Consider the random variables $(R_n)_{n\geq 1}$. Then, for   every $p> 2  $
$$\frac{1}{(n \log n)^p}\sum_{j=1}^n R_j I(R_j\leq t_j)\to 0, \qquad P-{\rm a.s.}$$
 where $t_n = 2^{m_n} (\log 2^{m_n})^\alpha$ for $0<\alpha<2$. 
\end{corollary}
\begin{proof}
 Take $f(x) = x \log^\alpha x$,   with $x\geq 1$ and   $0<\alpha<2$ and set  $R^\prime_j = R_j I(R_j\leq t_j)$ and $S_n^\prime = \sum_{j=1}^n R^\prime_j$. Let $p>2$ and take $q$ large enough so that  $2q(p-2)\geq1$ (which means $p \geq   2+\frac{1}{2q}$). Then 
\begin{align*}&
P\big(|S_n^\prime-a_n|\geq \varepsilon  a_n ^p\big)= P\big(|S_n^\prime-a_n|^{2q}\geq \varepsilon^{2q} a_n ^{2pq}\big) \\&\leq \frac{1}{\varepsilon^{2q} a_n ^{2pq}}E\big[(S_n^\prime-a_n)^{2q}\big]\leq \frac{2^{2q-1}}{\varepsilon^{2q} a_n ^{2pq}} E\big[(S_n^\prime)^{2q}\big] + \frac{2^{2q-1}}{\varepsilon^{2q} a_n ^{2q(p-1)}} \\&=   \frac{2^{2q-1}}{\varepsilon^{2q} a_n ^{2pq}} E   \left[ \left(\sum_{j=1}^n R^\prime_j\right)^{2q}\right]+ \frac{2^{2q-1}}{\varepsilon^{2q} a_n ^{2q(p-1)}}\leq  \frac{2^{2q-1}n^{2q-1}}{\varepsilon^{2q} a_n ^{2pq}}\left(\sum_{j=1}^nE\big[(R_j^\prime)^{2q}\big]\right)+ \frac{2^{2q-1}}{\varepsilon^{2q} a_n ^{2q(p-1)}}\\&\leq \frac{2^{2q-1}}{\varepsilon^{2q}}\left(\frac{1}{n^{2qp-2q +1}(\log n)^{2pq}}\sum_{j = 1}^n t^{2q}_j+ \frac{1}{(n\log n)^{2q(p-1)}}\right).
\end{align*} 
  The result follows by applying the Borel-Cantelli lemma since by Lemma \ref{lemma4}
	\[
	\sum_{n}\frac{1}{n^{2qp-2q +1}(\log n)^{2pq}}\sum_{j = 1}^n t^{2q}_j
	\]
	converges, while 
	\[
	\sum_{n}\frac{1}{(n\log n)^{2q(p-1)}}
	\]
	again converges since $ 2q(p-1)>2q(p-2)\geq 1$.
 \end{proof}

\medskip

 \noindent
We are ready for the proof of Theorem \ref{conv1}.

\bigskip
\noindent
\begin{proof}The proof is motivated by  the proof of Theorem 1 in \cite{M}.  In detail, take  $f(x) = x \log^\alpha x$,   for $x\geq 1$  with $\frac{1}{r}< \alpha < 1$   and notice that $f$ verifies the assumptions of   Lemma \ref{lemma1} (since $\alpha r >1$).  Recall the notation used before, i.e. 
$$t_n = f(2^{m_n}), \quad R^\prime_n =R_n I(R_n\leq t_n), \quad a_n = n\log n \quad\mbox{and}\quad S^\prime_n = \sum_{i=1}^{n}R^\prime_j.$$ 
Furthermore, for every $\varepsilon> 0$ put 
$$S_n(\varepsilon)=\sum_{j=1}^n R_j I(R_j\leq \varepsilon a_n).$$
Since
$$t_n = f(2^{m_n})\leq f (2^{\log_2 n})=f(n) =  n \log^\alpha n,$$
we have that
$$\lim_{n \to \infty} \frac{t_n}{a_n}=0.$$
 Thus,  for fixed  $\varepsilon >0$ we can  take $n$ sufficiently large such that $t_n < \varepsilon a_n$. Then
\begin{align*}&
|S_n(\varepsilon) - S^\prime_n|= \big|\sum_{j=1}^n R_j I(t_j <R_j\leq \varepsilon a_n)\big|\leq \varepsilon a_n N_{m_n}+ \sum_{k=1}^{m_n}\varepsilon a_{2^{m_n-k+1}}N_{m_n -k}\\
&\leq \varepsilon a_n \left(N_{m_n} + \sum_{k=1}^{m_n} \Big(\frac{1}{2}\Big)^{k-1}N_{m_n -k}\right) \leq \varepsilon  a_n^p \left(N_{m_n} + \sum_{k=1}^{m_n} \Big(\frac{1}{2}\Big)^{k-1}N_{m_n -k}\right)\\& \leq \varepsilon a_n^p N_{m_n}\left(1 + \sum_{k=1}^{m_n} \Big(\frac{1}{2}\Big)^{k-1}\right),
\end{align*}
where the third relation is due to the inequality $\frac{a_n}{a_{2n}} \leq \frac{1}{2}$.   Take any $s\geq r$; then, by Lemma      \ref{lemma5}, we obtain, ultimately   
$$ |S_n(\varepsilon) - S^\prime_n|\leq  \varepsilon a_n^{p} s \left(1 + \sum_{k=1}^{m_n} \Big(\frac{1}{2}\Big)^{k-1}  \right)\leq  3\varepsilon a_n^{p} s, \qquad P-{\rm a.s.} $$ 
 Moreover, 
\begin{align*}
|S_n(\varepsilon) - ^{(r)}S_n| = |\sum_{i=1}^{n}R_jI(R_j> \varepsilon a_n)- \sum_{k=2}^{r} M_n^{(k)}|\leq \sum_{i=1}^{n}R_jI(R_j> \varepsilon a_n)+ \sum_{k=2}^{r} M_n^{(k)}. 
\end{align*}
By \eqref{A3}, the first summand is finite for sufficiently large $n$ while Corollary \ref{corollary2} ensures that $\dfrac{1}{a_n}\sum_{k=2}^{r} M_n^{(k)} \to 0$. Then, as $n\to\infty$ 
\[
a_n^{-1}|S_n(\varepsilon) - ^{(r)}S_n|\to 0
\]
and so $|S_n(\varepsilon) - ^{(r)}S_n|\leq \varepsilon a_n$.  Finally, 
$$| ^{(r)}S_n  - S^\prime_n|\leq |S_n(\varepsilon) - S^\prime_n|+|S_n(\varepsilon) - ^{(r)}S _n|\leq \varepsilon a_n^p (3s+1)$$
and Corollary \ref{corollary1} gives the conclusion, by the arbitrariness of $\varepsilon.$ 
\end{proof}

 \medskip
 \noindent As a direct consequence of the Theorem obtained above we have the following asymptotic result for the particular case where $r =1$.
	\begin{theorem}\label{conv}\sl Consider the random variables $(R_n)_{n\geq 1}$ and assume that for the involved distribution functions $(F_n)_{n\geq 1}$ the following condition is satisfied:
		\begin{equation*}
		 \sup_{n \geq 1 }\limsup_{x\to 0}  \frac{F_n(x)}{x}<\infty.
		\end{equation*}Then, for every $p > 2$,
		\begin{equation}\label{claim}
		\frac{ S_n-M_n^{(1)}}{(n \log n)^p}\to 0, \qquad  P-{\rm a.s.} . 
		\end{equation}
	\end{theorem}
\begin{proof}
	 First, observe that for any $r \geq 2$,
	\[
	\frac{ S_n-M_n^{(1)}}{(n \log n)^p} = \dfrac{^{(r)}S_n}{(n\log n)^p} +\sum_{k=2}^{r} \dfrac{M_n^{k}}{(n\log n)^p}.
	\]
	The convergence of the latter expression is established by Theorem \ref{conv1} and Corollary \ref{corollary2}. 
\end{proof}
 
\begin{remark}
	\sl As mentioned at the beginning of the section, although the last result is weaker than the one proven in the previous section (observe that the convergence is to zero and not to a positive constant), it is obtained without imposing any conditions on the structure of the random variables $R_n$. Moreover, the involved distributions are not assumed to follow the same law and condition \eqref{as3F} is relaxed to condition \eqref{A1}.
\end{remark}

\end{document}